\documentclass[a4paper,11pt]{article}
\usepackage[latin1]{inputenc}
\usepackage{amsfonts}
\usepackage{amsmath}
\usepackage{amssymb}
\usepackage{amsthm}
\usepackage[T1]{fontenc}
\usepackage[dvips]{graphicx}
\usepackage{color}
\usepackage{epsfig}
\theoremstyle{plain}
\newtheorem{Thm}{Theorem}
\newtheorem{Lem}{Lemma}
\newtheorem{Def}{Definition}
\newtheorem*{Def*}{Definition}

\newtheorem*{Rem*}{\textsc{Remark}}
\newtheorem*{Lem*}{\textsc{Lemma}}
\newtheorem*{Cor*}{\textsc{Corollary}}
\newtheorem*{Con*}{\textsc{Conjecture}}
\newtheorem{Ass}{\textsc{Assumptions}}
\newcommand{\e}{\varepsilon}
\newcommand{\bb}[1]{\mathbb{ #1 }}
\newcommand{\av}[1]{\left| #1 \right|}

\title{Y-system and Deformed Thermodynamic Bethe Ansatz}
\author{Davide Masoero \thanks{E-mail address: masoero@sissa.it}\\ SISSA - Trieste}
\date{}
\begin{document}

\maketitle

\abstract{We introduce a new tool, the Deformed TBA (Deformed Thermodynamic Bethe Ansatz),
to analyze the monodromy problem of the cubic oscillator.
The Deformed TBA is a system of five coupled nonlinear integral equations, which in a particular case
reduces to the Zamolodchikov TBA equation
for the 3-state Potts model. Our method generalizes the Dorey-Tateo analysis of the (monomial) cubic oscillator.
We introduce a Y-system corresponding to the Deformed TBA and give it an elegant geometric interpretation.}

\section{Introduction}

This is the first in a series of papers that we dedicate to studying
the exact theory of the direct and inverse monodromy problem
for the cubic anharmonic oscillator 
\begin{equation}\label{eq:schr}
\frac{d^2\psi(\lambda)}{d\lambda^2}= V(\lambda;a,b) \psi(\lambda)\; , \quad
V(\lambda;a,b)=4 \lambda^3 - a \lambda - b \;.
\end{equation}
The main purpose of the present paper is to introduce a novel instrument of analysis,
that we call Deformed Thermodynamic Bethe Ansatz (Deformed TBA).

The monodromy problem of the cubic anharmonic oscillator is a fundamental and rather interesting problem in itself. 
Moreover, it is deeply interconnected with the study of singularities of solutions of the first Painlev\'e equations
(see author's papers \cite{piwkb}, \cite{piwkb2}).

Monodromy problems for anharmonic oscillators, and especially eigenvalue problems,
has been intensively studied since early '70 years with a wealth of different methods:
see for example the seminal papers  \cite{wu}, \cite{simon}, \cite{voros83}, \cite{bender98},
\cite{tateoTQ}, \cite{bazhanov00}, \cite{eremenko}.

Computations are often performed by means of the complex WKB methods (see for example
\cite{delabaere}, \cite{bender} and author's papers \cite{piwkb}, \cite{piwkb2})
or refined asymptotic expansions (see \cite{zinnjustin09} and references therein).
The first breaktrough towards an exact evaluation of
the monodromy problem is the work
of Dorey and Tateo \cite{tateoTQ}: they analyze anharmonic oscillators with a monomial potential
$\lambda^n-E$ ($n$ not necessarily 3)
via the Thermodynamic Bethe Ansatz and other nonlinear integral equations
(called sometimes Destri - de Vega equations).
Subsequently Bazhanov, Lukyanov and Zamolodchikov
generalized the Dorey-Tateo analysis to monomial potentials with a centrifugal term \cite{bazhanov00}.

Here we generalize Dorey and Tateo approach to the general cubic potential. In Theorem \ref{thm:dtba} below, we show
that the monodromy problem for the general cubic potential
is encoded in a \textit{nonlinear nonlocal Riemann-Hilbert problem}, which is equivalent
(at least for small value of the parameter $a$ in (\ref{eq:schr}))
to the following system of nonlinear integral equations that we call Deformed Thermodynamic Bethe Ansatz:

\begin{equation}\label{eq:dtba}
\chi_l(\sigma)=\int_{-\infty}^{+\infty}\varphi_l(\sigma-\sigma')\Lambda_l(\sigma') \,, \sigma,\sigma' \in \bb{R} \,,
\; l \in \bb{Z}_5 = \left\lbrace -2, \dots,2 \right\rbrace \; .
\end{equation}

Here

\begin{eqnarray*}
\Lambda_l(\vartheta)=\sum_{k \in \bb{Z}_5}e^{i\frac{2 l k \pi}{5}}L_k(\sigma) \; &,&
 L_k(\sigma)= \ln\left(1+e^{-\e_k(\sigma)}  \right) \, ,\\
 \e_k(\sigma) \!=\frac{1}{5}\sum_{l \in \bb{Z}_5}e^{-i\frac{2 l k\pi}{5}}\chi_l(\sigma) \!\! \!&\!+\!&\! 
 \frac{\sqrt{\frac{\pi}{3}}\Gamma(1/3)}{2^{\frac{5}{3}} \Gamma(11/6)} e^{\sigma}\!+ \!
 a \!\frac{\sqrt{3 \pi}\Gamma(2/3)}{4^{\frac{2}{3}} \Gamma(1/6)} e^{\frac{\sigma}{5}-i\frac{2 k \pi}{5}} \,, \\
\varphi_0(\sigma) =\frac{\sqrt{3}}{\pi}\frac{2 \cosh(2\sigma)}{1+ 2\cosh(2\sigma)} &,&
\varphi_1(\sigma) =-\frac{\sqrt{3}}{\pi}
\frac{e^{-\frac{9}{5}\sigma}}{1+ 2\cosh(2\sigma)} \\ 
\varphi_2(\sigma) =-\frac{\sqrt{3}}{\pi}
\frac{e^{-\frac{3}{5}\sigma}}{1+ 2\cosh(2\sigma)} \!\!\!\!&\!,\!&\!\!\!\!
\varphi_{-1}(\sigma) = \varphi_1(-\sigma) \,, \varphi_{-2}(\sigma) = \varphi_2(-\sigma)\,.
\end{eqnarray*}

Here the pseudo-energy $\e_k$ is related to
the logarithm of the k-th Stokes multiplier (see Section 3 below).

For $a=0$, see equation (\ref{eq:tba}) below, equations (\ref{eq:dtba}) reduce to the
Thermodynamic Bethe Ansatz, introduced by Zamolodchikov \cite{zamo1} to describe the thermodynamics of
the 3-state Potts model and Lee-Yang model.

The paper is organized as follows.
In Section 2 we give a geometric construction of the \textit{space of monodromy data}.
We realize the Stokes multipliers of the cubic oscillators
as natural coordinates on the quotient $W_5/PSL(2,\bb{C})$,
where $W_5$ is a dense  open subset of $\left( \bb{P}^1 \right)^5$ (the Cartesian product of five copies of
$\bb{P}^1$).
Section 3 is devoted to the construction of the deformed Y-system.
In Section 4 we derive the Deformed Thermodynamic Bethe Ansatz.
For convenience of the reader, we explain the basic theory of cubic oscillators (Stokes sectors, Stokes multipliers,
subdominant solutions, etc ...) in the Appendix.

\paragraph{Acknowledgments}
I am indebted to my advisor Prof. B. Dubrovin who constantly gave me
suggestions and advice. This work began in December 2009 during a visit at the "Dipartimento di Fisica Teorica"
of Turin University. I thank Dr R. Tateo for the kind hospitality, for useful discussions and encouragement.
I also thank Prof. Y. Takei and Dr M. Mazzocco for continuous interest in my work.
This work is partially supported by the Italian Ministry of University and Research
(MIUR) grant PRIN 2008 "Geometric methods in the theory of nonlinear waves and their applications".

\section{The space of monodromy data}
Usually monodromy data of equation (\ref{eq:schr}) are expressed in terms of Stokes multipliers, which
are defined by means of a special set of solutions of the equation. In this section,
following Nevanlinna \cite{nevanlinna} and author's paper \cite{piwkb},
we study the monodromy data from a geometric (hence invariant) viewpoint.
Here and for the rest of the paper $\bb{Z}_5=\left\lbrace-2,\dots,2\right\rbrace$.
\begin{Def}
Let $\left\lbrace \varphi,\chi \right\rbrace$ be a basis of solution of (\ref{eq:schr}).

We  call \begin{equation}\label{def:wk}
 w_k(\varphi,\chi)=\lim_{\substack{\lambda\ \to \infty \\
\left|\arg{\lambda}-\frac{2 \pi k}{5}\right| < \frac{\pi}{5} -\varepsilon}}\frac{\varphi(\lambda)}{\chi(\lambda)} \in \mathbb{C}
\cup \infty \, , \; k \in \mathbb{Z}_5 \; .
\end{equation}
the k-th asymptotic value.
\end{Def}

We collect the main properties of the asymptotic values in the following
\begin{Lem}\label{lem:wk}

\begin{itemize}
 \item[(i)]  Let $\varphi'= a \, \varphi + b \, \chi$ and $\chi'=c \, \varphi + d \, \chi'$,
$\left(\begin{matrix}
a & b \\ 
c  & d
\end{matrix} \right) \in Gl(2,\bb{C})$. Then
\begin{equation}\label{eq:moebius}
w_k(\varphi',\chi')= \frac{a \, w_k(\varphi,\chi) + b}{c \, w_k(\varphi,\chi)+d} \; .
\end{equation}
\item[(ii)] $w_{k-1}(\varphi,\chi)=w_{k+1}(\varphi,\chi)$ \textit{iff} $\sigma_k(a,b)=0$ . Here $\sigma_k(a,b)$ is the k-th
Stokes multipliers.
\item[(iii)]$w_{k+1}(\varphi,\chi) \neq w_{k}(\varphi,\chi)$

\end{itemize}

\begin{proof}
See \cite{piwkb}.
\end{proof}
\end{Lem}

\begin{Def}\label{def:w5}
We define $$W_5=\left\lbrace (z_{-2},z_{-1},z_{0},z_{1},z_{2}), z_k  \in \bb{C}\cup \infty,
z_k \neq z_{k \pm 1}\right\rbrace \; .$$
\end{Def}
The group of automorphism of the Riemann sphere, called M\"obius group or $PSL(2,\bb{C})$,
has the following natural free action
on $W_5$: let $T = \left(
\begin{matrix}
a & b \\ 
c  & d
\end{matrix} \right) \in PSL(2,\bb{C})$ then
$$T(z_{-2}, \dots, z_2)=(\frac{a z_{-2} +b}{c z_{-2} +d},\dots, \frac{a z_{2} +b}{c z_{2} +d}) \; .$$

After formula (\ref{def:wk}) and Lemma \ref{lem:wk}(iii) every basis of solution of (\ref{eq:schr}) determines a point in $W_5$. After the transformation law (\ref{eq:moebius}), the Schr\"odinger equation (\ref{eq:schr}) determines an orbit
of the $PSL(2,\bb{C})$ action. Hence, we define the \textit{space of monodromy data} as follows.

\begin{Def}\label{def:v5}
We call \textit{space of monodromy data} the space of the orbits of the $PSL(2,\bb{C})$ action
and denote it $V_5$:
\begin{equation*}
V_5=W_5/PSL(2,\bb{C}) \;.
\end{equation*}
\end{Def}

\begin{Rem*}
The author proved in \cite{piwkb} that the space of monodromy data $V_5$ is the moduli space of
solutions of the first Painlev\'e equation. $V_5$ is a smooth manifold and ${\cal{M}}_{0,5} \subset V_5 \subset \overline{{\cal{M}}_{0,5}}$.
\end{Rem*}

In the rest of the section we construct natural coordinates of $V_5$ and interpret them
in the spirit of cubic oscillator theory.
Define
\begin{eqnarray}\nonumber
R_k: W_5 &\to& \mathbb{C} \, , \; k \in \bb{Z}_5  \; ,\\ \label{eq:Rk}
R_k(z_{-2},\dots,z_{2}) &=& \left(z_{1+k},z_{-2+k}; z_{-1+k},z_{2+k} \right) \;,
\end{eqnarray}
where $(a,b;c,d)=\frac{(a-c)(b-d)}{(a-d)(b-c)}$ is the cross ratio of four point on the sphere.

We collect the main properties of functions $R_k$ in the following Lemma, whose easy proof is left to the reader.

\begin{Lem}\label{lem:Rk}
 \begin{itemize}
   \item[(i)]
   \begin{eqnarray}\nonumber
  R_k(z_{-2},\dots,z_2) &\neq& \infty \,, \; \forall (z_{-2},\dots,z_2) \in W_5 \;, \\ \label{eq:Rkvalues}
  R_k(z_{-2},\dots,z_2)&=&0 \mbox{ iff } z_{k-1}=z_{k+1} \;, \\ \nonumber
  R_k(z_{-2},\dots,z_2)&=&0 \mbox{ iff } z_{k-1}=z_{k+2} \mbox{ or } z_{k+1}=z_{k-2} \; .
           \end{eqnarray}
   \item[(ii)] The functions $R_k$ are invariant under the $PSL(2,\bb{C})$ action. Hence they 
   are well defined on $V_5$: with a small abuse of notation we denote $R_k$ also the functions defined
  on  $V_5$.
  \item[(iii)] They satisfy the following set of quadratic relation \footnote{almost identical to the relations
  satisfied by the Stokes multipliers (\ref{eq:multipliers})}
  \begin{equation}\label{eq:quad}
R_{k-2}R_{k+2}=1 - R_k \,, \; \forall k \in \bb{Z}_5 \, .
 \end{equation}
 \item[(iv)]The pair $R_k, R_{k+1}$ is a coordinate system of $V_5$ on the open subset $R_{k-2} \neq 0$.
 The pair of coordinate systems $(R_k,R_{k+1})$ and $(R_{k+2},R_{k-2})$ form an atlas of $V_5$.
 \end{itemize}

\end{Lem}

\begin{Def}\label{def:map}
We call any cubic polynomial of the form $V(\lambda;a,b)=4\lambda^3 -a\lambda -b$ a cubic potential.
The above formula identifies the space of cubic potentials with $\mathbb{C}_2 \backepsilon (a,b)$.
Through the asymptotic values (\ref{def:wk}), we define the monodromy map 
\begin{equation*}
{\cal{T}}: \mathbb{C}^2 \to V_5  \; .
\end{equation*}
${\cal{T}}(a,b)$ is the monodromy data of equation (\ref{eq:schr}).
\end{Def}
Nevanlinna \cite{nevanlinna} proved that the monodromy map $\cal{T}$ is surjective.

With a small abuse of notation we define
\begin{equation}\label{eq:potRk}
R_k(a,b)= R_k \circ {\cal{T}}(a,b) \; .
\end{equation}

The reader can verify that the Stokes multipliers (defined precisely in the Appendix) are (modulo multiplicative constant)
the functions $R_k$:
$$
\sigma_k(a,b)=i R_k(a,b) \; .
$$
We have thus realized the Stokes multipliers geometrically as natural coordinates of the monodromy map from
the space of cubic oscillators to the quotient $W_5/PSL(2,\bb{C})$.

\begin{Rem*}
The same construction presented here holds for anharmonic oscillators with polynomial potentials
of any degree. If $n$ is the degree, we denote $V_{n+2}$ the space of monodromy data and $R_k^{(n+2)}$ the
natural functions (defined by formula \ref{eq:Rk}). For example $V_3$ is one point and $V_4$ is $\bb{C}^*$.
The functions $R_k^{(n+2)}$ satisfy a system of algebraic relations similar to (\ref{eq:potRk}), which will
be studied in a subsequent paper.
The reader should notice that only in the cubic case the functions $R_k^{(n+2)}$ coincide with the Stokes multipliers.
\end{Rem*}

\section{Y-system}
Here we introduce the Y-system (\ref{eq:Ysystem}), which is a fundamental step
in the derivation of the Deformed TBA.

We begin with an observation, probably due to Sibuya \cite{sibuya75}:

\begin{Lem}\label{lem:0monodromy}
Let $\omega=e^{i\frac{2 \pi}{5}}$ and $R_k$ be defined as in (\ref{eq:potRk}). Then
\begin{equation}\label{eq:0monodoromy}
R_k(\omega^{-1}a,\omega b)=R_{k-2}(a,b) \; .
\end{equation}

\begin{proof}
Denote $\varphi(\lambda;a,b)$ a solution of (\ref{eq:schr}) whose Cauchy data
do not depend on $a,b$. It is an entire function of three complex variables
with some remarkable properties.
For any $k \in \mathbb{Z}_5$ $\varphi(\omega^k\lambda;\omega^{2k}a,\omega^{3k}b)$ satisfies the same Schr\"odinger equation (\ref{eq:schr}) as $\varphi(\lambda;a,b)$.
Fix $\varphi(\lambda;a,b),\chi(\lambda;a,b)$ linearly independent solutions and define
\begin{eqnarray*}
w_k(a,b)&=&w_k(\varphi(\lambda;a,b),\chi(\lambda;a,b)) \; , \\
\bar{w}_k(a,b)&=&w_l(\varphi(\omega^l\lambda;\omega^{2l}a,\omega^{3l}b),\chi(\omega^l\lambda;\omega^{2l}a,\omega^{3l}b) \; .
\end{eqnarray*}
Then $w_k(\omega^{2l}a,\omega^{3l}b)=\bar{w}_{k-l}(a,b)$. Choose $l=2$ to obtain the thesis.
\end{proof}
\end{Lem}

Due to equations (\ref{eq:quad}) and relations (\ref{eq:0monodoromy}), the holomorphic functions $R_k(a,b)$
satisfies the following system of functional equations, first studied by Sibuya \cite{sibuya75}

\begin{equation}\label{eq:Sibuya}
R_k(\omega^{-1}a,\omega b)R_k(\omega a, \omega^{-1}b)= 1 -R_k(a,b) \, , \; \forall k \in \bb{Z}_5.
\end{equation}

We have collected all the elements to introduce the important \textit{Y-functions} and \textit{Y-system}.

We fix $a \in \bb{C}$ and define
\begin{equation}\label{eq:Yk}
Y_k(\vartheta)=-R_0(\omega^{-k}a,e^{\frac{6}{5}\vartheta})  \, ,\; k \in \bb{Z}_5 \; .
\end{equation}

Sibuya's equation (\ref{eq:Sibuya}) is equivalent to
the following system of functional equations, that we call \textit{Deformed Y-system}:

\begin{equation}\label{eq:Ysystem}
Y_{k-1}(\vartheta-i\frac{\pi}{3})Y_{k+1}(\vartheta+i\frac{\pi}{3})=1+Y_{k}(\vartheta) \; .
\end{equation}

\begin{Rem*}
If $a=0$, $Y_k=Y_0 , \forall k$ and the system (\ref{eq:Ysystem}) reduces to just one equation, called
$Y$-system, which was introduced by Zamolodchikov \cite{zamo2} in relation with the Lee-Yang and 3-state Potts models.
Dorey and Tateo \cite{tateoTQ} studied the Zamolodchikov $Y$-system in relation with the
Schr\"odinger equation with potential $V(\lambda;0,b)=4 \lambda^3 -b$.
\end{Rem*}

\subsection{Analytic Properties of $Y_k$}
In the following theorem we summarize the analytic properties of the Y-functions.

\begin{Thm}\label{thm:Y}
\begin{itemize}
\item[(i)] For any $a \in \bb{C}$ and $k \in \bb{Z}_5$, $Y_k$ is analytic and $i \frac{5\pi}{3}$ periodic.
If $a$ is real then $\overline{Y_{k}(\overline{\vartheta})}=Y_{-k}(\vartheta)$,
where $\overline{~^{~}}$ stands for complex conjugation.
\item[(ii)]For any $a \in \bb{C}$ and $k \in \bb{Z}_5$,
\begin{eqnarray}\nonumber
\av{\frac{Y_k(\vartheta)}{\tilde{Y}_k(\vartheta)} - 1} &=& O (e^{-\mbox{Re}\vartheta}), \mbox{ as } Re{\vartheta} \to +\infty \mbox{ and } \av{Im{\vartheta}} < \frac{\pi}{2} \,, \\ \label{eq:dominantY} 
\tilde{Y}_k(\vartheta)&=& \exp\left( A e^{\vartheta} + B a e^{\frac{\vartheta}{5}-i\frac{2 k \pi}{5}}\right) \; .
\end{eqnarray}
Here $A=\frac{\sqrt{\frac{\pi}{3}}\Gamma(1/3)}{2^{\frac{5}{3}} \Gamma(11/6)}$ and
$B =\frac{\sqrt{3 \pi}\Gamma(2/3)}{4^{\frac{2}{3}} \Gamma(1/6)}$.
\item[(iii)]For any $a \in \bb{C}$ and any $K \in \bb{R}$, $Y_k(\vartheta)$ is bounded on $Re{\vartheta} \leq K$.
If $a=0$, $\lim_{\vartheta \to -\infty}Y_k(\vartheta)=\frac{1+ \sqrt{5}}{2}$.
\item[(iv)] If $e^{i\frac{2 k \pi}{5}} a$ is real non negative then
$Y_k(\vartheta)=0$ implies $Im{\vartheta}=\pm \frac{5 \pi}{6}$.
If $a=0$ then $Y_k(\vartheta)=-1$ implies $\vartheta=\pm i\frac{\pi}{2}$
\item[(v)] Fix $\e > 0$. If $a$ is small enough, then for any $k \in \bb{Z}_5$, $Y_k(\vartheta) \neq 0,-1$
for any $\vartheta \in \av{\mbox{Im} \vartheta} \leq \frac{\pi}{2} -\e$.
\end{itemize}
\begin{proof}
\begin{itemize}
 \item[(i)]Trivial.
 \item[(ii)]These "WKB-like" estimates can be found in \cite{piwkb} Section 4 or in \cite{sibuya75}.
 \item[(iii)]The boundedness follows directly from the fact that $R_k(a,b)$ is analytic in $b=0$.
 If $a=b=0$, then for symmetry reasons one can choose $\varphi, \chi$ such that $w_k=e^{i\frac{2 k \pi}{5}}$.
 This implies the thesis.
 \item[(iv)]The statement is equivalent to Theorem \ref{thm:pt} in the Appendix.
 \item[(v)]Since $Y_k$ depends analytically from the parameter $a$, it follows from (iv).
\end{itemize}

\end{proof}

\end{Thm}

\section{Deformed TBA}
This section is devoted to the derivation
of the Deformed Thermodynamic Bethe Ansatz equations (\ref{eq:dtba}) \footnote{The reader who wants to repeat all the computations below, should remember that $\sum_{l=0}^{n-1}e^{i\frac{2 \pi l}{n}}=0$.}.

In what follows we always make the following
\begin{Ass}\label{assump}
We assume that there exists an $\e>0$ such that
\begin{itemize}
 \item[(i)]every branch  of $\ln{Y_k}$  is holomorphic on $\av{Im \vartheta} \leq \frac{\pi}{3} + \e$,
 and bounded for $\vartheta \to - \infty$. And
 \item[(ii)]every branch of $\ln{(1+\frac{1}{Y_k})}$ is holomorphic on $\av{Im \vartheta} \leq + \e$,
 and bounded for $\vartheta \to - \infty$.
\end{itemize}
\end{Ass}

From Theorem \ref{thm:Y}(iii, v) we know that the assumptions are valid
if $a$ is small enough.

We define the following bounded analytic functions on \textit{the physical strip}
$Im{\vartheta}\leq\frac{\pi}{3}$
\begin{eqnarray} \label{eq:pseudoen}
\e_k(\vartheta) &=& \ln{Y_k(\vartheta)}, \\ \nonumber
\delta_k(\vartheta)&=&\e_k(\vartheta)- \frac{\sqrt{\frac{\pi}{3}}\Gamma(1/3)}{2^{\frac{5}{3}} \Gamma(11/6)}
 e^{\vartheta} + a  \frac{\sqrt{3 \pi}\Gamma(2/3)}{4^{\frac{2}{3}} \Gamma(1/6)}
 e^{\frac{\vartheta}{5}-i\frac{2 k \pi}{5}},\\ \nonumber
L_k(\vartheta)&=&\ln(1+e^{-\e_k(\vartheta)}) \; .
\end{eqnarray}
Here the branches of logarithms are fixed by requiring
\begin{equation}\label{eq:asymptotic}
\lim_{\sigma \to +\infty}\delta_k(\sigma + i \tau)= \lim_{\sigma \to +\infty}L_k(\sigma + i \tau)=  0 \, ,
\; \forall \av{\tau} <\frac{\pi}{2} \; .
\end{equation}
We remark that by Theorem \ref{thm:Y}(ii), this choice is always possible.
We denote $\e_k$ pseudo-energies in analogy with the undeformed TBA.

Due to the Y-system (\ref{eq:Ysystem}), the functions $\delta_k$ satisfy the following
\textit{nonlinear nonlocal Riemann-Hilbert problem}
 \begin{equation}\label{eq:r-h}
\delta_{k-1}(\vartheta -i\frac{\pi}{3})+ 
\delta_{k+1}(\vartheta +i\frac{\pi}{3})-\delta_k(\vartheta)=L_k(\vartheta) \, , \; \av{Im\vartheta} \leq \e \; .
 \end{equation}
Here the boundary conditions are given by asymptotics (\ref{eq:asymptotic}).

The system  (\ref{eq:r-h}) is $\bb{Z}_5$ invariant. Hence we diagonalize its linear part (the left hand side)
by taking its \textit{discrete Fourier transform} (also called \textit{Wannier} transform):
\begin{eqnarray}\label{eq:chi}
\chi_l(\vartheta) = \sum_{k \in \bb{Z}_5}e^{i\frac{2 kl\pi}{5}}\delta_k(\vartheta) \!\! &,&
\delta_k(\vartheta) = \frac{1}{5}\sum_{l \in \bb{Z}_5}e^{-i\frac{2 kl\pi}{5}}\chi_l(\vartheta) \, ,\\ \nonumber
\Lambda_l(\vartheta)=\sum_{k \in \bb{Z}_5}e^{i\frac{2 l\pi}{5}}L_k(\vartheta) \!\! &,&
L_k(\vartheta) = \frac{1}{5}\sum_{l \in \bb{Z}_5}e^{-i\frac{2 kl\pi}{5}}\Lambda_l(\vartheta) \;.
\end{eqnarray}

The above defined functions satisfy the following system

\begin{equation}\label{eq:diagonal}
e^{-i\frac{2 l \pi}{5}}\chi_l(\vartheta +i\frac{\pi}{3})+ e^{i\frac{2l \pi}{5}}
\chi_l(\vartheta -i\frac{\pi}{3})-\chi_l(\vartheta)=\Lambda_l(\vartheta).
\end{equation}

The system of functional equation (\ref{eq:diagonal}), may be rewritten in the convenient form of a system of coupled integral equations (\ref{eq:dtba}).
\begin{Thm}\label{thm:dtba}
If $a$ is small enough, the functions $\chi_l$ satisfy the Deformed Thermodynamic Bethe Ansatz
\begin{equation*}
\chi_l(\sigma)=\int_{-\infty}^{+\infty}\varphi_l(\sigma-\sigma')\Lambda_l(\sigma') \,, \sigma,\sigma' \in \bb{R} \;.
\; \; \quad \quad \quad (\ref{eq:dtba}) 
\end{equation*}

Here $\Lambda_l$ are defined as in (\ref{eq:pseudoen},\ref{eq:chi}) and

\begin{eqnarray}\nonumber
\varphi_0(\sigma) &=&\frac{\sqrt{3}}{\pi}\frac{2\cosh(\sigma)}{1+2\cosh(2\sigma)} \\ \nonumber
\varphi_1(\sigma) &=&-\frac{\sqrt{3}}{\pi}
\frac{e^{-\frac{9}{5}\sigma}}{1+2\cosh(2\sigma)} \\ \label{eq:kernels}
\varphi_2(\sigma) &=&-\frac{\sqrt{3}}{\pi}
\frac{e^{-\frac{3}{5}\sigma}}{1+2\cosh(2\sigma)} \\ \nonumber
\varphi_{-1}(\sigma) &=&-\frac{\sqrt{3}}{\pi}
\frac{e^{\frac{9}{5}\sigma}}{1+2\cosh(2\sigma)} \\ \nonumber
\varphi_{-2}(\sigma) &=&-\frac{\sqrt{3}}{\pi}
\frac{e^{\frac{3}{5}\sigma}}{1+2\cosh(2\sigma)} \,.
\end{eqnarray}

\begin{proof}
If $a$ is small enough we know that the Assumptions \ref{assump} are valid.
Hence the thesis follows from system (\ref{eq:diagonal}) and the technical Lemma \ref{lem:r-h} below.
\end{proof}

\end{Thm}

\begin{Lem}\label{lem:r-h}
Let  $f: \av{Im \vartheta} \leq \e \to \mathbb{C} $ be a bounded analytic function. Then for any $l \in \bb{Z}_5$ there exists a unique function $F$ analytic and bounded  on $\av{Im \vartheta }\leq\frac{\pi}{3} + \e$,
such that $$e^{-i\frac{2 \pi l}{5}}F(\vartheta+i \frac{\pi}{3})+ e^{i\frac{2 \pi l}{5}} F(\vartheta-i \frac{\pi}{3})- F(\vartheta)=f(\vartheta) \,,\forall \av{Im \vartheta} \leq \e \; .$$
Moreover, $F$ is expressed through the following integral transform
\begin{equation}\label{eq:transform}
F(\vartheta + i \tau)=\int_{-\infty}^{+\infty}\varphi_l(\vartheta+ i \tau -\vartheta')f(\vartheta')d\sigma' \, , \; \forall \av{Im \vartheta} \leq \e, \av{\tau} \leq \frac{\pi}{3}\; \, ,
\end{equation}
provided $\av{Im(\vartheta+ i \tau -\vartheta')} < \frac{\pi}{3}$ and the integration path belongs to the strip $\av{Im \vartheta'}\leq \e$. Here $\varphi_l$ is defined by formula (\ref{eq:kernels}).

\begin{proof}
Uniqueness: let $F_1, F_2$ be bounded solution  of the functional equation
$$e^{-i\frac{2 \pi l}{5}}F_j(\vartheta+i \frac{\pi}{3})+ e^{i\frac{2 \pi l}{5}} F_j(\vartheta-i \frac{\pi}{3})- F_j(\vartheta)=f(\sigma) \, , \; j=1,2 , \av{Im\vartheta}\leq \e \;.$$
Their difference $G=F_1-F_2$ satisfies $$e^{-i\frac{2 \pi l}{5}}G(\vartheta+i \frac{\pi}{3})+ e^{i\frac{2 \pi l}{5}} G(\vartheta-i \frac{\pi}{3})-G(\vartheta)=0 \, .$$ Hence $G$ extends to an entire, $2 i \pi$ periodic, bounded function.
Therefore $G$ is a constant and the only constant satisfying the functional relation is zero.

Existence: one notices that if $\theta \neq \pm i n \frac{\pi}{3},\, n \in \bb{Z}$ then
$$
e^{-i\frac{2 \pi l}{5}}\varphi_l(\theta+i \frac{\pi}{3})+ e^{i\frac{2 \pi l}{5}} \varphi_l(\theta-i \frac{\pi}{3})
- \varphi_l(\theta) =0 \, , \; \forall l \in \mathbb{Z}_5 \;.$$ Then a rather standard computation shows that the function $F$ defined through formula (\ref{eq:transform}) satisfies all the desired properties.
\end{proof}

\end{Lem}

\begin{Rem*}
Once the system of integral equations (\ref{eq:dtba}) is solved for $\sigma \in \bb{R}$, one can use
the same set of integral equations as explicit formulas to extend the functions
$\chi_l(\theta)$ on $\av{Im \theta} \leq \frac{\pi}{3}$.
Then one can use the Y-system (\ref{eq:Ysystem}) to extend the Y functions on the entire fundamental strip
$\av{Im \vartheta} \leq \frac{5 \pi}{6}$.
\end{Rem*}

\begin{Rem*}
While the Y-system equations do not depend on the parameter $a$ (the coefficient of the linear term
of the potential $4\lambda^3- a\lambda -b$),
on the contrary the Deformed TBA equations depend on it since it enters explicitely
into the definition of functions $\Lambda_l$.
\end{Rem*}

\subsection{The case $a=0$}

If $a=0$ then $\delta_k=\delta_0, \, L_k=L_0\forall k$.
Therefore, $\delta_0$ satisfy the single functional equation
\begin{equation}\label{eq:a=0r-h}
\delta_{0}(\vartheta -i\frac{\pi}{3})+ 
\delta_{0}(\vartheta +i\frac{\pi}{3})-\delta_0(\vartheta)=L_0(\vartheta) \, , \; \av{Im\vartheta} \leq \e \; .
\end{equation}

Similar reasoning as in Theorem \ref{thm:dtba} shows that $\delta_0$ satisfies
the following nonlinear integral equation (as it was firstly discovered by Dorey and Tateo \cite{tateoTQ})
\begin{equation}\label{eq:tba}
\delta_0(\sigma)=\int_{-\infty}^{+\infty}\varphi_0(\sigma-\sigma')\ln{(1+\exp{(-(\delta_0(\sigma')+A e^{\sigma'})})} \,, \sigma,\sigma' \in \bb{R} \;.
\end{equation}
Here $\varphi_0$ is defined as in formula (\ref{eq:kernels}) and
$A=\frac{\sqrt{\frac{\pi}{3}}\Gamma(1/3)}{2^{\frac{5}{3}} \Gamma(11/6)}$.

Equation (\ref{eq:tba}) is called
Thermodynamic Bethe Ansatz and was introduced by Zamolodchikov \cite{zamo1} to describe the Thermodynamics of
the 3-state Potts and Lee-Yang models. 

\section{Concluding Remarks}
We have given a geometric construction of the \textit{space of monodromy data} of cubic oscillators
and we have used such construction to analyze the direct monodromy problem.

We have shown that for small value of the deformation parameter (the coefficient of the linear term of the potential)
the direct monodromy problem
defines a nonlinear nonlocal Riemann-Hilbert problem which is equivalent
to a system of nonlinear equations, that we call Deformed Thermodynamic Bethe Ansatz (Deformed TBA).
Our approach can be easily generalized to anharmonic oscillators of any order. We will give the details
in a subsequent publication (see also the Remark at the end of Section 2).

As it was said in the introduction, this is just the first of a series of paper
that we dedicate to the Deformed TBA. Our rather ambitious purpose is to use the Deformed TBA to
construct analytical and numerical tools to effectively solve the monodromy problem of
the general cubic oscillator. We plan to pursue our study of the Deformed TBA equations
in different directions to achieve such a goal.

In particular, together with A. Moro and R. Tateo we are implementing an algorithm to solve numerically the equations
and to understand for which values of the parameter $a$ the Assumptions \ref{assump} fail.
Consequently, we will modify the Riemann-Hilbert problem for taking into account the multivaluedness of
the functions $\Lambda_l, \chi_l$.

At the same time, we are going to study analytically the convergence of successive approximation
schemes to solve the Deformed TBA.

As it is known from author previous works \cite{piwkb},\cite{piwkb2} (see also \cite{takei}),
cubic oscillators are deeply interconnected with the first Painlev\'e equation and in particular with
the distribution of poles of its solutions. As an application of the theory developed here,
we are going to investigate the poles of special solution of first Painlev\'e equation. In particular we will study
the tritronqu\'ee solution,
since it is relevant to describe the singular limit of the focusing Nonlinear Schr\"odinger equation \cite{dubrovin},
\cite{bertola10}.

\section{Appendix}

The reader expert in anharmonic oscillators theory will skip this Appendix; for her, it will be enough to know
that we denote $\sigma_k(a,b)$ the $k-th$ Stokes multipliers of equation (\ref{eq:schr}).
Here we review briefly the standard way, i.e. by means of Stokes multipliers,
of introducing the monodromy problem
for equation (\ref{eq:schr}). All the statements of this section are proved in
Appendix A of author's paper \cite{piwkb} and in Sibuya's book \cite{sibuya75} .

\begin{Lem}\label{lem:wkb}
Fix $k \in \mathbb{Z}_5 = \left\lbrace -2, \dots , 2 \right\rbrace$,
define the branch of $\lambda^{\frac{1}{2}}$ by requiring
$$\lim_{\substack{\lambda \to \infty \\ \arg{\lambda}=
\frac{2 \pi k}{5}}} {\rm Re}\lambda^{\frac{5}{2}} = + \infty \,$$
while choose the branch of $\lambda^{\frac{1}{4}}$ globally on the complex plane minus the negative
real axis such that it is positive on the positive real axis.
Then there exists a unique solution $y_k(\lambda;a,b)$ of equation (\ref{eq:schr})
such that 
\begin{equation}\label{eq:asym}
\lim_{\substack{ \lambda \to \infty \\ \av{\arg{\lambda} - \frac{2 \pi k}{5}} < \frac{3 \pi}{5} -\e}}
\frac{y_k(\lambda;a,b)}{\lambda^{-\frac{3}{4}} e^{-\frac{4}{5} \lambda^{\frac{5}{2}} +
\frac{a}{2}\lambda^{\frac{1}{2}}}} \to 1 , \; \forall \e >0 \, .
\end{equation}

\end{Lem}

\begin{Def*}
We denote $\left\lbrace \lambda \in \bb{C}, \av{\arg \lambda - \frac{2 k \pi}{5} } < \frac{\pi}{5} \right\rbrace$
the k-th Stokes sector.
We denote $y_k$, defined in Lemma above, the k-th subdominant solution or the solution subdominant in the k-th sector.
\end{Def*}

From the asymptotics (\ref{eq:asym}), it follows that $y_k$ and $y_{k+1}$
are linearly independent and the following equations hold true:

\begin{eqnarray}\nonumber
 y_{k-1}(\lambda;a,b) &=& y_{k+1}(\lambda;a,b)+\sigma_k(a,b)y_k(\lambda;a,b) \, ,\\ \label{eq:multipliers} \\ \nonumber
 - i \sigma_{k+3}   &=& 1+\sigma_k(a,b)\sigma_{k+1}(a,b)\, , \; \forall k \in \bb{Z}_5 \; .
\end{eqnarray}

\begin{Def*}
The entire functions $\sigma_k(a,b)$ are called Stokes multipliers.
The quintuplet of Stokes multipliers $\sigma_k(a,b), k \in \bb{Z}_5$ is called the monodromy data
of equation (\ref{eq:schr}).
\end{Def*}

\subsection{PT-symmetric anharmonic oscillator}
The eigenvalues problem $\sigma_k=0, k \in \bb{Z}_5$ is PT symmetric if $\omega^k a \in \bb{R}$
(the study of PT symmetric oscillators began in the seminal paper \cite{bender98}).
Dorey, Dunning, Tateo \cite{tateo} and Shin \cite{shin} proved the following result about PT
symmetric cubic oscillator.

\begin{Thm}\label{thm:pt}
Fix $k \in \bb{Z}_5$. Suppose $\sigma_k(a,b)=0$ and $\omega^{2k} a$ is real and \textit{positive}.
Then $\omega^{3k}b$ is real and negative.
\end{Thm}


\bibliographystyle{alpha}
\newcommand{\etalchar}[1]{$^{#1}$}

\end{document}